\documentclass[a4paper, article, oneside, UKenglish]{amsart}

\usepackage[colorinlistoftodos,prependcaption,textsize=tiny]{todonotes}

\usepackage[utf8]{inputenx} 
\usepackage[T1]{fontenc}    

\usepackage{lmodern}           
\usepackage[scaled]{beramono}  

\vfuzz \hfuzz
\usepackage{amsmath}
\usepackage{amssymb}     
\usepackage{amsthm}    
\usepackage{thmtools}  
\usepackage{mathtools} 
\usepackage{mathrsfs}  

\usepackage{graphicx}  
\usepackage{babel}     
\usepackage{csquotes}  
\usepackage{textcomp}  
\usepackage{listings}  
\usepackage{varioref}
\usepackage{hyperref}
\usepackage[nameinlink, capitalize, noabbrev]{cleveref}

\declaretheorem[style = plain, numberwithin = section]{theorem}
\declaretheorem[style = plain,      sibling = theorem]{corollary}
\declaretheorem[style = plain,      sibling = theorem]{lemma}
\declaretheorem[style = plain,      sibling = theorem]{proposition}
\declaretheorem[style = definition, sibling = theorem]{definition}
\declaretheorem[style = definition, sibling = theorem]{example}

\declaretheorem[style = remark,    numbered = yes]{remark}

\DeclareMathOperator{\di}{\delta}

\DeclareMathOperator{\ord}{ord}

\newcommand{\N}{\mathbb{N}}   
\newcommand{\C}{\mathbb{C}}   
\renewcommand{\P}{\mathbb{P}} 

\newcommand{\ol}{\overline}
\newcommand{\PP}{\mathbb{P}^1\times\mathbb{P}^1}
\newcommand{\V}{V}
\newcommand{\noi}{\noindent}

\renewcommand{\a}{\alpha}
\renewcommand{\b}{\beta}
\renewcommand{\c}{c}

\def\di{\partial}

\def\Ø{\emptyset}

\crefname{listing}{Program}{Programs}
\Crefname{listing}{Program}{Programs}

\usepackage{booktabs}
\usepackage{subfigure}

\begin{document}
\title[Weierstrass points on algebraic curves in $\PP$]{Special Weierstrass points\\on algebraic curves in $\PP$}
\author{Paul Aleksander Maugesten} 
\address{Department of Mathematics\\  University of Oslo\\ P.O. Box 1053 Blindern\\ NO-0316 Oslo\\ NORWAY} 
\email{paulamau@math.uio.no} 
\author{Torgunn Karoline Moe} 
\address{The Science Library\\ University of Oslo\\ P.O. Box 1063 Blindern\\ NO-0316 Oslo\\ NORWAY} 
\email{torgunnk@math.uio.no}
\maketitle
\begin{abstract}
    \noindent
In this paper we demonstrate that the notion of inflection points and extactic points on plane algebraic curves can be suitably transferred to curves in $\PP$. 

More precisely, we describe osculating curves and study Weierstrass points of algebraic curves in the surface $\PP$ with respect to certain linear systems. In particular, we study points where a fiber of $\PP$ is tangent, and points with a hyperosculating $(1,1)$-curve. In the first case we find Hessian-like curves that intersect the curve in these points, and in the second case we find a local criteria. Moreover, we provide Plücker-like formulas for the number of smooth Weierstrass points on a curve. In the special case of rational curves, we use suitable Wronskians to compute these points and their respective Weierstrass weights.
    
\end{abstract}


\tableofcontents

\section{Introduction}
The surface $\PP$ is a well known surface with simple properties; a first example of a projective surface after the projective plane $\P^2$. However, curves in this surface have not been subject to the same thorough exploration as plane curves. 

Our first encounter with curves in $\PP$ was in the Ph.D. thesis of the second author \cite{Moe13}, where the main focus is cuspidal curves, i.e. curves with only unibranched singularities. Osculating properties of points on curves in $\PP$ is mentioned there on a side note, and in this paper we pick up some of these previously unpublished results and take a new look at them from the point of view of our recent article \cite{MM17}. In \cite{MM17} we study plane algebraic curves, sextactic points and 2-Hessian curves, including a Plücker-like counting formula, and it turns out that similar features are possible to describe for curves in $\PP$ with quite classical methods. In this article we develop these results for the complete linear systems of $(0,1)$-, $(1,0)$- and $(1,1)$-curves on curves in $\PP$. 

The article has the following structure. In \cref{sec:back} we dwell on some important background results that will be used throughout the article. In \cref{sec:fibertangent} we derive Plücker-like formulas and corresponding Hessian curves for the first two cases. In \cref{sec:oneone} we similarly find a formula for the number of smooth $(1,1)$-Weierstrass points on cuspidal curves, and we calculate a local criteria to determine if a point on a curve is a $(1,1)$-Weierstrass point. The results follow from methods by \cite{BG97} and \cite{Not99}. In \cref{sec:rational} we study rational curves and obtain similar results as in \cref{sec:oneone} using suitable Wronski determinants. Lastly, in \cref{sec:applications} we apply our results to an example and push a bit further.

\section{Background}\label{sec:back}
In this section we recall the basic properties of $\PP$ and curves in this surface.

\subsection{Curves in $\PP$}
With $\C[x_0,x_1;y_0,y_1]$ the coordinate ring (Cox ring) of $\PP$, a curve $C$ of type $(a,b)$ in ${\PP}$ is given as the zero set of a bigraded, bihomogeneous polynomial \[F(x_0,x_1,y_0,y_1) \in \C[x_0,x_1;y_0,y_1]_{(a,b)}, \quad a,b \geq 0.\] 

With $C$ an integral algebraic curve in ${\PP}$, a point $p$ on $C$ is singular if all the partial derivatives vanish, smooth otherwise. In particular, a unibranched singularity is referred to as a cusp. Let $m_p$ denote the multiplicity of a point $p$, and $\delta_p$ its delta invariant.

A first classification of curves in $\PP$ can be done via the geometric genus $g$. 
\begin{proposition}[Genus formula]
	\label{for:genus}
	The genus $g$ of a curve $C$ of type $(a,b)$ in $\PP$ is given by
	\begin{equation*}
	g=(a-1)(b-1)-\sum \delta_p,
	\end{equation*}
	where the sum is taken over all points $p \in C$.
	\begin{proof}
		See \cite[p.~393]{Har77} and \cite[Corollary 3.1.4]{Moe13}.
	\end{proof}
	
\end{proposition}
Note that integral curves of type $(a,1)$ and $(1,b)$ are non-singular.

\subsection{Osculating curves in $\PP$}
Now we transfer the notion of osculating and hyperosculating curves from $\P^2$ to $\PP$. 

For a fixed pair $(\alpha,\beta) \in \mathbb{N}^2$, consider curves, not necessarily integral, of type $(\alpha,\beta)$ in $\PP$ and their intersections with $C$. Let \[r(\alpha,\beta)=(\alpha+1)(\beta+1)-1,\] i.e. the projective dimension of $H^0(\PP,\mathcal{O}_{\PP}(\alpha,\beta))$. 

\begin{definition}
A curve $T^{(\alpha,\beta)}_p$ of type $(\alpha,\beta)$ is called an osculating curve to $C$ at a smooth point $p$ if it intersects $C$ at $p$ with local intersection multiplicity \[(T^{(\alpha,\beta)}_p\cdot C)_p \geq r(\alpha,\beta).\]
When the inequality is strict, the curve $T^{(\a,\b)}_p$ is said to be hyperosculating, and $p$ an $(\a,\b)$-Weierstrass point.
\end{definition}

Note that a singular point of $C$ is always considered an $(\a,\b)$-Weierstrass point, see \cite{BG97}.

Lifting to a global perspective, we say that a curve with the property that it intersects $C$ in its $(\a,\b)$-Weierstrass points, is an $(\a,\b)$-Hessian to $C$, denoted by $H_{(\a,\b)}$. 

\subsection{Euler relations}\label{euler}
The classical Euler relation in $\P^2$ has its equivalents in $\PP$. These relations will be used in proofs throughout the article.
\begin{lemma}\label{euler1}
	Let $C=V(F)$ be a curve of type $(a,b)$ in $\PP$. For any point $p \in C$ with coordinates $(a_0:a_1;b_0:b_1)$ the following relations hold:
	\begin{align*}
	a_0F_{x_0}(p)+a_1F_{x_1}(p)&=aF(p)=0,\\
	b_0F_{y_0}(p)+b_1F_{y_1}(p)&=bF(p)=0.
	\end{align*}
\begin{proof}
The result follows by evaluating at $p$ the following identities for a polynomial $F \in \C[x_0,x_1,y_0,y_1]_{(a,b)}.$
\begin{align*}
x_0F_{x_0}+x_1F_{x_1}&=aF,\\
y_0F_{y_0}+y_1F_{y_1}&=bF.
\end{align*}
\end{proof}
\end{lemma}
By differentiating, combining and evaluating the identities in the proof of \cref{euler1}, we get a number of new identities and, by applying them to points on $C$, new relations. We list the main results in the form of identites, and note that the list is not in any way exhaustive.
\begin{lemma}\label{euler2} 
For a polynomial $F \in \C[x_0,x_1,y_0,y_1]_{(a,b)}$, the following identities hold:	
\begin{align*}
x_0 F_{x_0x_0}+x_1 F_{x_0x_1}&=(a-1)F_{x_0}\\
x_0 F_{x_0x_1}+x_1 F_{x_1x_1}&=(a-1)F_{x_1}\\
x_0^2 F_{x_0x_0}+2x_0x_1 F_{x_0x_1}+x_1^2 F_{x_1x_1}&=(a-1)aF\\
y_0 F_{y_0y_0}+y_1 F_{y_0y_1}&=(b-1)F_{y_0}\\
y_0 F_{y_0y_1}+y_1 F_{y_1y_1}&=(b-1)F_{y_1}\\
y_0^2 F_{y_0y_0}+2y_0y_1 F_{y_0y_1}+y_1^2F_{y_1y_1}&=(b-1)bF\\
x_0 F_{x_0y_0}+x_1 F_{x_1y_0}&=aF_{y_0}\\
x_0 F_{x_0y_1}+x_1 F_{x_1y_1}&=aF_{y_1}\\
y_0 F_{x_0y_0}+y_1 F_{x_0y_1}&=bF_{x_0}\\
y_0 F_{x_1y_0}+y_1 F_{x_1y_1}&=bF_{x_1}\\
x_0y_0 F_{x_0y_0}+x_1y_0 F_{x_1y_0}+x_0y_1 F_{x_0y_1}+x_1y_1 F_{x_1y_1}&=abF.\\
\end{align*}
\end{lemma}

\section{$(1,0)$- and $(0,1)$-Weierstrass points}\label{sec:fibertangent}
In this section we explore the complete linear systems of $(1,0)$- and $(0,1)$-curves on curves of type $(a,b)$ in $\PP$, given by \begin{align*}
\gamma_{x_0}x_0+\gamma_{x_1}x_1&=0, \quad \gamma_{x_0},\gamma_{x_1} \in \C,\\
\gamma_{y_0}y_0+\gamma_{y_1}y_1&=0,  \quad \gamma_{y_0},\gamma_{y_1} \in \C.
\end{align*} These curves will frequently be referred to as fibers, reflecting the fact that $\PP$ is a doubly ruled surface over $\P^1$.

\subsection{Osculating $(1,0)$- and $(0,1)$-curves}
The first result is trivial, and reveals that the fibers through a point are the osculating curves in this case. 
\begin{proposition}[Osculating $(1,0)$- and $(0,1)$-curves]\label{osc1001curve}
The osculating $(1,0)$-curve to a curve $C$ at a point $p=(a_0:a_1;b_0:b_1)$ is given by \[a_1x_0-a_0x_1=0, \text{ or equivalently, } F_{x_0}(p)x_0+F_{x_1}(p)x_1=0,\] and similarly, the osculating $(0,1)$-curve is given by \[b_1y_0-b_0y_1=0, \text{ or equivalently, } F_{y_0}(p)y_0+F_{y_1}(p)y_1=0.\] 
\begin{proof}
	Since $r(1,0)=1=r(0,1)$, the respective unique fibers passing through a point $p$ are the osculating $(1,0)$- and $(0,1)$-curves to $C$ at $p$. By \cref{euler1}, \[a_0F_{x_0}(p)+a_1F_{x_1}(p)=0 \text{ and } b_0F_{y_0}(p)+b_1F_{y_1}(p)=0,\] and the result follows.
\end{proof}
\end{proposition}

\subsection{$(1,0)$- and $(0,1)$-Weierstrass points}
We now turn to points where a $(1,0)$- or $(0,1)$-curve is hyperosculating.
\begin{proposition}
	A point $p$ on a curve $C=V(F)$ is a $(1,0)$-Weierstrass point if and only if $F_{y_0}(p)=0=F_{y_1}(p)$. Similarly, $p$ is a $(0,1)$-Weierstrass point if and only if $F_{x_0}(p)=0=F_{x_1}(p)$.
	\begin{proof}
	We move $p$ to $(0:1;0:1)$ by a linear transformation, and designate $x_0$ and $y_0$ as local coordinates. Since $F_{y_0}(p)=0$, \[F(x_0,1,y_0,1)=  F_{x_0}(p)x_0+H(2),\] where $H(2)$ denotes higher order terms, and consequently the fiber $T_p^{(1,0)}=V(x_0)$ through $p$ intersects $C$ at $p$ with multiplicity $(T^{(1,0)}_p \cdot C)_p\geq 2>r(1,0)=1$.
	\end{proof}
\end{proposition}

\begin{remark}
A singular point is both a $(1,0)$- and $(0,1)$-Weierstrass point.
\end{remark}

\begin{remark}
Abusing definitions, we sometimes refer to the property of a fiber considered as a curve in an affine neighbourhood of $p$, and say that a fiber is tangent to $C$ at a point $p$ if the fiber $T_p$ is a tangent line in the classical sense, i.e. $(T_p \cdot C)_p>m_p$, where $m_p$ is the multiplicity of $p$.	
\end{remark}

\begin{corollary}\label{cor:1001points}	
A point $p=(a_0:a_1;b_0:b_1)$ on a curve $C$ of type $(a,b)$, $a,b\geq2$, is a $(1,0)$-Weierstrass point if and only if \[F_{y_0y_0}(p)F_{y_1y_1}(p)-F_{y_0y_1}^2(p)=0.\] Similarly, $p$ is a $(0,1)$-Weierstrass point if and only if \[F_{x_0x_0}(p)F_{x_1x_1}(p)-F_{x_0x_1}^2(p)=0.\]  	
\begin{proof}
We prove this for $(0,1)$-Weierstrass points, the other case is symmetrical.

Assume first that $a_0,a_1 \neq 0$. Then by the Euler relations, $F_{x_0}(p)=0$ if and only if $F_{x_1}(p)=0$, so $p$ is a $(0,1)$-Weierstrass point if and only if $F_{x_0}(p)F_{x_1}(p)=0$. Exploiting the Euler identities in \cref{euler2}, we find that:
{\small{
\begin{align*}
(a-1)^{2}F_{x_0}F_{x_1}&=(x_0F_{x_0x_0}+x_1F_{x_0x_1})(x_0F_{x_0x_1}+x_1F_{x_1x_1})\\
&=(x_0^2F_{x_0x_0}F_{x_0x_1}+x_0x_1F_{x_0x_1}^2+x_1^2F_{x_0x_1}F_{x_1x_1}+x_0x_1F_{x_0x_0}F_{x_1x_1})\\
&=F_{x_0x_1}(x_0^2F_{x_0x_0}+2x_0x_1F_{x_0x_1}+x_1^2F_{x_1x_1})+x_0x_1(F_{x_0x_0}F_{x_1x_1}-F_{x_0x_1}^2)\\
&=x_0x_1(F_{x_0x_0}F_{x_1x_1}-F_{x_0x_1}^2),
\end{align*}}}so in this case $F_{x_0}(p)F_{x_1}(p)=0$ if and only if $F_{x_0x_0}(p)F_{x_1x_1}(p)-F_{x_0x_1}^2(p)=0$.

Next, assume that $a_0=0$ and $F_{x_0}(p)=0=F_{x_1}(p)$. Consider the two differentiated Euler identities 
\begin{align*}
(a-1)F_{x_0}&=x_0F_{x_0x_0}+x_1F_{x_0x_1},\\ 
(a-1)F_{x_1}&=x_0F_{x_1x_0}+x_1F_{x_1x_1}.
\end{align*}
 Evaluating these identities at $p$ gives $F_{x_0x_1}(p)=0=F_{x_1x_1}(p)$. In particular, \[F_{x_0x_0}(p)F_{x_1x_1}(p)-F_{x_0x_1}^2(p)=0.\] Conversely, assume that $a_0=0$ and $F_{x_0x_0}(p)F_{x_1x_1}(p)-F_{x_0x_1}^2(p)=0$. Since $a_0=0$, $a_1 \neq 0$, so the relation $0=aF(p)=a_1F_{x_1}(p)$, and $F_{x_1}(p)=0$. Similarly, $0=(a-1)F_{x_1}(p)=a_1F_{x_1x_1}(p)$, hence $F_{x_1x_1}(p)=0$. Combining this with the second assumption, we have $F_{x_0x_1}^2(p)=0$, which reduces to $F_{x_0x_1}(p)=0$. Thus, $(a-1)F_{x_0}(p)=a_1F_{x_0x_1}(p)=0$, and since $a \geq 2$, we conclude that $F_{x_0}(p)=0$. A parallell argument holds in the case that $a_1=0$.
\end{proof}
\end{corollary}

With the above result in mind, candidates for $(1,0)$- and $(0,1)$-Hessian curves follow suit.
\begin{theorem}[The $(1,0)$- and $(0,1)$-Hessian]
	Let $C$ be a curve of type $(a,b)$, with $a,b\geq2$. Then the $(1,0)$-Hessian curve to $C$ is given by \[F_{y_0y_0}F_{y_1y_1}-F_{y_0y_1}^2=0.\]
Similarly, the $(0,1)$-Hessian curve to $C$ is given by \[F_{x_0x_0}F_{x_1x_1}-F_{x_0x_1}^2=0.\]  
\begin{proof}
By \cref{cor:1001points}, these curves intersect $C$ only at the corresponding $(1,0)$- and $(0,1)$-Weierstrass points.	
\end{proof}
\end{theorem}

Mixing the two directions, it is possible to find the $(1,0)$- and $(0,1)$-Weierstrass points with a combined criteria.
\begin{theorem}
A point $p$ on $C$ is either a $(1,0)$- or a $(0,1)$-Weierstrass point if and only if \[F_{x_0y_0}(p)F_{x_1y_1}(p)-F_{x_0y_1}(p)F_{x_1y_0}(p)=0.\]

\begin{proof}
If $p=(a_0:a_1;b_0:b_1)$ is a $(1,0)$-Weierstrass point, then
\begin{align*}
aF_{y_0}(p)&=0=a_0F_{x_0y_0}(p)+a_1F_{x_1y_0}(p),\\ 
aF_{y_1}(p)&=0=a_0F_{x_0y_1}(p)+a_1F_{x_1y_1}(p).
\end{align*} 
Assume first that $a_0,a_1 \neq 0$. Then sorting these relations and multiplying each side yields \[a_0a_1F_{x_0y_0}(p)F_{x_1y_1}(p)=a_0a_1F_{x_0y_1}(p)F_{x_1y_0}(p),\] and the relation in the theorem holds. If $a_0=0$ and $a_1 \neq 0$, then the equations lead to \[F_{x_1y_0}(p)=0=F_{x_1y_1}(p),\] and again the relation holds. A similar argument can be used if $a_1=0$. The case of $(0,1)$-Weierstrass points follows by symmetry. 

Conversely, assume that \[F_{x_0y_0}(p)F_{x_1y_1}(p)-F_{x_0y_1}(p)F_{x_1y_0}(p)=0\] and $p$ is not a $(1,0)$-Weierstrass point. We need to show that it is a $(0,1)$-Weierstrass point. If $a_0,a_1,b_0,b_1 \neq 0$, then {\small{
\begin{align*}
0&=\left(\tfrac{a}{a_0}F_{y_0}(p)-\tfrac{a_1}{a_0}F_{x_1y_0}(p)\right)\left(\tfrac{a}{a_1}F_{y_1}(p)-\tfrac{a_0}{a_1}F_{x_0y_1}(p)\right)-F_{x_0y_1}(p)F_{x_1y_0}(p)\\
&=\tfrac{a^2}{a_0a_1}F_{y_0}(p)F_{y_1}(p)-\tfrac{a}{a_0}F_{x_1y_0}(p)F_{y_1}(p)-\tfrac{a}{a_1}F_{x_0y_1}(p)F_{y_0}(p)\\
&=\tfrac{a}{a_0a_1}\bigl(aF_{y_0}(p)F_{y_1}(p)-{a_1}F_{x_1y_0}(p)F_{y_1}(p)-{a_0}F_{x_0y_1}(p)F_{y_0}(p)\bigr).\\
\end{align*}}} In this case, the expression above can be rewritten in two ways by using that at least one of $F_{y_0}(p)$ and $F_{y_1}(p)$ is non-zero, say the latter.

First, observe that
\begin{align*}
{a_0}F_{x_0y_1}(p)F_{y_0}(p)&=\bigl(aF_{y_0}(p)-{a_1}F_{x_1y_0}(p)\bigr)F_{y_1}(p)\\
{a_0}F_{x_0y_1}(p)F_{y_0}(p)&={a_0}F_{x_0y_0}(p)F_{y_1}(p)\\
-\tfrac{b_1}{b_0}F_{x_0y_1}(p)F_{y_1}(p)&=F_{x_0y_0}(p)F_{y_1}(p)\\
{b_0}F_{x_0y_0}(p)+{b_1}F_{x_0y_1}(p)&=0\\
bF_{x_0}(p)&=0\\
F_{x_0}(p)&=0.
\end{align*}
And similarly,
\begin{align*}
{a_1}F_{x_1y_0}(p)F_{y_1}(p)&=\bigl(aF_{y_1}(p)-{a_0}F_{x_0y_1}(p)\bigr)F_{y_0}(p)
\end{align*} implies $F_{x_1}(p)=0$. Hence, the point $p$ is a $(0,1)$-Weierstrass point.

Secondly, assume that $a_0=0$ and $a_1 \neq 0$, and both $b_0,b_1 \neq 0$. In this case, $F_{x_1}(p)=0$, and subsequently $b_0F_{x_1y_0}(p)=-b_1F_{x_1y_1}(p)$. Hence, if the relation holds,
\begin{align*}
0=&F_{x_0y_0}(p)F_{x_1y_1}(p)-F_{x_0y_1}(p)F_{x_1y_0}(p)\\
=&-\tfrac{1}{b_1}F_{x_1y_0}(p)\bigl(b_0F_{x_0y_0}(p)+b_1F_{x_0y_1}(p)\bigr)\\
=&-\tfrac{b-1}{b_1}F_{x_1y_0}(p)F_{x_0}(p),
\end{align*}
If $F_{x_0}(p)=0$, then $p$ is a $(0,1)$-Weierstrass point, contrary to assumption. If $F_{x_0}(p)\neq0$, we must have $F_{x_1y_0}(p)= 0$. Therefore \[aF_{y_0}(p)=a_0 F_{x_0y_0}(p)+a_1F_{x_1y_0}(p)=0,\] and since $a \geq 2$,  $F_{y_0}(p)=0$. Moreover, since $b_1\neq 0$ and $b_0F_{y_0}(p)+b_1F_{y_1}(p)=0$, it follows that $F_{y_1}(p)=0$, thus $p$ is a $(1,0)$-Weierstrass point. 

Thirdly, assume that both $a_0,a_1 \neq 0$, and $b_0=0$ and $b_1 \neq 0$. Then $F_{y_1}(p)=0$ and subsequently, \[a_0F_{x_0y_1}(p)+a_1F_{x_1y_1}(p)=0.\] As above, the relation gives 
\begin{align*}
0=&-\tfrac{1}{a_1}F_{x_0y_1}(p)\bigl(a_0F_{x_0y_0}(p)+a_1F_{x_1y_0}(p)\bigr)\\
=&-\tfrac{a-1}{a_1}F_{x_0y_1}(p)F_{y_0}(p).
\end{align*}
If $F_{y_0}(p)=0$, then we are done. If $F_{x_0y_1}(p)=0$, then \[bF_{x_0}(p)=b_0 F_{x_0y_0}(p)+b_1F_{x_0y_1}(p)=0,\] and since $b \geq 2$,  $F_{x_0}(p)=0$. Moreover, since $a_1\neq 0$ and $a_0F_{x_0}(p)+a_1F_{x_1}(p)=0$, it follows that $F_{x_1}(p)=0$, thus $p$ is a $(0,1)$-Weierstrass point. 

Lastly, if $p=(0:1;0:1)$, then $F_{x_1}(p)=0=F_{y_1}(p)$ and furthermore {${F_{x_1y_1}(p)=0}$}. Then either $F_{x_0y_1}(p)=0$ and thus $F_{x_0}(p)=0$, or $F_{x_1y_0}(p)=0$ and thus $F_{y_0}(p)=0$.
\end{proof}
\end{theorem}

\begin{corollary}\label{cor:1001Hessian}
Let $C$ be a curve of type $(a,b)$, with $a,b\geq2$. Then  \[F_{x_0y_0}F_{x_1y_1}-F_{x_0y_1}F_{x_1y_0}=0\] defines a curve $H_{(1,0)\cup(0,1)}$ that intersects $C$ in its $(1,0)$- and $(1,0)$-Weierstrass points.
\end{corollary}

\subsection{The number of $(1,0)$- and $(0,1)$-Weierstrass points}
We now find a formula for the number of smooth $(1,0)$- and $(0,1)$-Weierstrass points on a curve.
\begin{theorem}\label{eq:formula1001} Let $C$ be a curve of type $(a,b)$ in $\PP$. For each point $p$ on $C$, let $(C,p_j)$ denote the branches of $C$ at $p$. Let $l_{p_j}$ denote the intersection multiplicity of the branch $(C,p_j)$ and the osculating $(0,1)$- or $(1,0)$-curve in the respective cases. 
	
		Then the number $W_{(1,0)}$ of $(1,0)$-Weierstrass points on $C$ (counted with multiplicity) is
	\begin{align*}
	W_{(1,0)}&=2(b+g-1)-\sum (l_{p_j}-1)\\
	&=2a(b-1) - 2\sum_{p \in C} \delta_p -\sum (l_{p_j}-1),
	\end{align*}
	where the sum $\sum (l_{p_j}-1)$ is taken over all branches of all singular points of $C$.
	
	Similarly, the number of $(0,1)$-Weierstrass points $W_{(0,1)}$ on $C$ is
	\begin{align*}
	W_{(0,1)}&=2(a+g-1)-\sum(l_{p_j}-1)\\
	&=2(a-1)b- 2\sum_{p \in C} \delta_p -\sum(l_{p_j}-1).
	\end{align*}
	
	\begin{proof}
		The complete linear system $Q$ of $(1,0)$-curves on $C$ has dimension $r=r(1,0)=1$, and $\deg Q=b$. By \cite[Proposition~3.4]{BG97}, the sum of the Weierstrass weights $w_p(Q)$ of the points in $C$ equals
		\begin{align*}
		\sum_{p\in C} w_p(Q) = 2(b + g - 1).
		\end{align*}
		
		By \cref{for:genus}, the right hand side can be rewritten as  
		\begin{align*}
		\sum_{p\in C} w_p(Q) &= 2\bigl(b + (a-1)(b-1)-\sum \delta_p - 1\bigr)\\
		&=2a(b-1)-2\sum \delta_p.
		\end{align*}
		
		For the left hand side, first assume that $p$ is unibranched. Then $w_p(Q)$ can be computed as 
		\begin{equation*}
		\label{eq:wp-used}
		w_p(Q) = \sum_{i=0}^1 (h_i - i),
		\end{equation*}
		where $h_i$ are distinct integers with the property that $h_i=(C \cdot C_i)_p$ for curves $C_0,C_1$ of type $(1,0)$, see \cite{Not99}.
		
		Trivially, $h_0=0$, and $h_1$ equals the intersection multiplicity $l_{p_j}$ of the branch and the unique $(1,0)$-curve through $p$. For singular points, summing up over all branches of $C$ gives the total $w_p(Q)$.
		
		Lastly, observe that for a smooth point, $w_p(Q) = l_p-1$, which equals its multiplicity as a $(1,0)$-Weierstrass point, and summing over all such points add up to $W_{(1,0)}$.
	\end{proof}
\end{theorem}

\begin{remark}
	Note that this result also follows from Hurwitz's theorem, see \cite[Theorem 3.1.6]{Moe13}.
	\end{remark}

\begin{remark}\label{conj:1001}
	The expressions in \cref{eq:formula1001} can be rewritten to the following form,  	
	\begin{align*}
2a(b-1)&=\sum_{p \in C} \bigl(2\delta_p+w_p(1,0)\bigr),\\
2(a-1)b&=\sum_{p \in C} \bigl(2\delta_p+w_p(0,1)\bigr).
	\end{align*} 
	Summing the two equations leads to
	\begin{equation}\label{eq:1001}
4ab-2a-2b=\sum_{p \in C} (4\delta_p+w_p(0,1)+w_p(1,0))
	\end{equation}
	The left hand side can be recognized as the intersection of $C$ (of type $(a,b)$) with a curve of type $(2(a-1),2(b-1))$, i.e. the curve $H_{(1,0)\cup(0,1)}$ of \cref{cor:1001Hessian}. The right hand side sums over all points in $C$, and we conjecture that the contribution from each point equals the local intersection multiplicity $(H_{(1,0)\cup(0,1)}\cdot C)_p$. In that case, \cref{eq:1001} reflects Bézout's theorem in $\PP$.
\end{remark}

\section{$(1,1)$-Weierstrass points}\label{sec:oneone}
In this section we consider $(1,1)$-curves in $\PP$, and develop similar results as in the previous section.

\subsection{$(1,1)$-curves in $\PP$}
A general $(1,1)$-curve in $\PP$ can be given by \[\gamma_{00}x_0y_0+\gamma_{01}x_0y_1+\gamma_{10}x_1y_0+\gamma_{11}x_1y_1=0, \quad \gamma_{ij} \in \C.\]

We first show that there is a unique $(1,1)$-curve through three generic points in $\PP$.
\begin{proposition}
	Let $p$, $q$ and $r$ be three distinct points in $\PP$ not all on the same $(1,0)$- or $(0,1)$-curve. The three points define a unique  $(1,1)$-curve $Q_{pqr}$. Moreover, if no pair of points are on the same $(1,0)$- or $(0,1)$-curve, then the $(1,1)$-curve is irreducible.
	
	In either case, if the points have coordinates 
	\begin{eqnarray*}
		p & = & (p_{x_0}:p_{x_1}; p_{y_0}:p_{y_1}),\\
		q & = & (q_{x_0}:q_{x_1}; q_{y_0}:q_{y_1}),\\
		r & = & (r_{x_0}:r_{x_1}; r_{y_0}:r_{y_1}),
	\end{eqnarray*}
	then the defining polynomial of $Q_{pqr}$ is given by
	\[G=\det \begin{bmatrix}
	x_0y_0 & x_0y_1 & x_1y_0 & x_1y_1 \\
	p_{x_0}p_{y_0} & p_{x_0}p_{y_1}  & p_{x_1}p_{y_0}  & p_{x_1}p_{y_1} \\
	q_{x_0}q_{y_0} & q_{x_0}q_{y_1}  & q_{x_1}q_{y_0}  & q_{x_1}q_{y_1} \\
	r_{x_0}r_{y_0} & r_{x_0}r_{y_1}  & r_{x_1}r_{y_0}  & r_{x_1}r_{y_1} 
	\end{bmatrix}.\]
	
	\label{11kurverlemma}
\begin{proof}
	There are two cases to consider. First, if any pair of points, say $p$ and $q$ are on the same $(1,0)$-curve, then the product of the defining polynomials of this curve and the unique $(0,1)$-curve through $r$ gives a unique, reducible $(1,1)$-curve passing through the three points. 

Second, if $p$, $q$ and $r$ are generic, then the points can be moved by a change of coordinates to
		\begin{eqnarray*}
			p & = & (1:0; 1:0),\\
			q & = & (0:1; 0:1),\\
			r & = & (1:1; 1:1).
		\end{eqnarray*}
		Let $$G=\gamma_{00}x_0y_0+\gamma_{01}x_0y_1+\gamma_{10}x_1y_0+\gamma_{11}x_1y_1$$ be the defining polynomial of a $(1,1)$-curve containing $p$, $q$ and $r$. Then $G(p)=0$ implies $\gamma_{00}=0$, and $G(q)=0$ implies $\gamma_{11}=0$. Moreover, $G(r)=0$ implies that $\gamma_{01}=-\gamma_{10}$. Thus, the coefficients of the defining polynomial of the $(1,1)$-curve containing $p$, $q$ and $r$ is uniquely determined up to multiplication with $c \in \C^*$. The curve must be irreducible, since if it were reducible, then, contrary to assumption, two of the points must have been on the same $(1,0)$- or $(0,1)$-curve. 

	The polynomial
	$$G=\det \begin{bmatrix}
	x_0y_0 & x_0y_1 & x_1y_0 & x_1y_1 \\
	p_{x_0}p_{y_0} & p_{x_0}p_{y_1}  & p_{x_1}p_{y_0}  & p_{x_1}p_{y_1} \\
	q_{x_0}q_{y_0} & q_{x_0}q_{y_1}  & q_{x_1}q_{y_0}  & q_{x_1}q_{y_1} \\
	r_{x_0}r_{y_0} & r_{x_0}r_{y_1}  & r_{x_1}r_{y_0}  & r_{x_1}r_{y_1} 
	\end{bmatrix}$$ 
	is the defining polynomial of a $(1,1)$-curve containing $p$, $q$ and $r$, and by uniqueness it is the defining polynomial of $Q_{pqr}$.
\end{proof}
\end{proposition}

\subsection{Osculating $(1,1)$-curves}
We now turn to osculating $(1,1)$-curves, and in this case $r=r(1,1)=3$.

\begin{theorem}\label{thm:osc11}
	Let $C$ be a curve of type $(a,b)$ in $\PP$, with $a+b \geq 3$. For every smooth point $p=(a_0:a_1;b_0:b_1)$ on a curve $C$ of type $(a,b)$ in $\PP$ there exists a unique $(1,1)$-curve $T^{(1,1)}_p$ such that $(T_p^{(1,1)} \cdot C)_p\geq 3$. 
	
	If a fiber is tangent to $C$ at $p$, then $T_p^{(1,1)}$ is the union of the two fibers through the point.
	
	If $a_i, b_j \neq 0$ and no fiber is tangent to $C$ at $p$, the defining polynomial of $T_p^{(1,1)}$ is given by \[\sum_{i,j=0}^1\Bigl(F_{x_iy_j}(p)-\frac{F_{x_ix_i}(p)F_{y_j}(p)}{2F_{x_i}(p)}-\frac{F_{y_jy_j}(p)F_{x_i}(p)}{2F_{y_j}(p)} \Bigr)x_iy_j.\]
	
	If one coordinate is zero, say $p=(0:1;b_0:b_1)$, where $b_0,b_1 \neq 0$, and no fiber is tangent to $C$ at $p$, then the defining polynomial of $T_p^{(1,1)}$ is given by
\[\sum_{j=0}^1\Bigl(F_{x_0y_j}(p)-\frac{F_{x_0x_0}(p)F_{y_j}(p)}{2F_{x_0}(p)}-\frac{F_{y_jy_j}(p)F_{x_0}(p)}{2F_{y_j}(p)} \Bigr)x_0y_j+F_{y_0}(p)x_1y_0+F_{y_1}(p)x_1y_1.\] Symmetrical relations hold for other combinations of vanishing coordinates.

In particular, if $p=(0:1;0:1)$ and no fiber is tangent to $C$ at $p$, then the defining polynomial of $T_p^{(1,1)}$ is given by
\[\Bigl(F_{x_0y_0}(p)-\frac{F_{x_0x_0}(p)F_{y_0}(p)}{2F_{x_0}(p)}-\frac{F_{y_0y_0}(p)F_{x_0}(p)}{2F_{y_0}(p)} \Bigr)x_0y_0+F_{x_0}(p)x_0y_1+F_{y_0}(p)x_1y_0.\]

\begin{proof}
	Assume first that $p=(0:1;0:1)$. A general $(1,1)$-curve \[\gamma_{00}x_0y_0+\gamma_{01}x_0y_1+\gamma_{10}x_1y_0+\gamma_{11}x_1y_1=0, \quad \gamma_{ij} \in \C,\] passes through $p$ if $\gamma_{11}=0$. The intersection at $p$ is non-transversal if the tangent lines in an affine neighbourhood coincide, which is fulfilled if $\gamma_{01}=F_{x_0}(p)$ and $\gamma_{10}=F_{y_0}(p)$.  
	
	To determine $\gamma_{00}$, we compute the intersection multiplicity at $p$ by means of local polynomials of $F$ and $T=\gamma_{00}x_0y_0+F_{x_0}(p)x_0y_1+F_{y_0}(p)x_1y_0$. Let $f(x_0,y_0)=F(x_0,1,y_0,1)$ and  $t(x_0,y_0)=T(x_0,1,y_0,1)$. 
		
	For two constants $k_1$ and $k_2$, put $P=f-t-k_1y_0t-k_2x_0t$. Then
		\begin{align*}
		P=H(3)&+\bigl(\tfrac{1}{2}F_{x_0x_0}(p)-k_2F_{x_0}(p)\bigr)x_0^2\\
		&+ \bigl(F_{x_0y_0}(p)-k_1F_{x_0}(p)-k_2F_{y_0}(p)-\gamma_{00}\bigr)x_0y_0\\
		&+\bigl(\tfrac{1}{2}F_{y_0y_0}(p)-k_1F_{y_0}(p)\bigr)y_0^2.
		\end{align*}
		
		The three coefficients in the above expression vanish whenever 
		\begin{align*}
		k_1&=\frac{F_{y_0y_0}(p)}{2F_{y_0}(p)},\\
		k_2&=\frac{F_{x_0x_0}(p)}{2F_{x_0}(p)},\\
		\gamma_{00}&=F_{x_0y_0}(p)-\frac{F_{y_0y_0}(p)F_{x_0}(p)}{2F_{y_0}(p)}-\frac{F_{x_0x_0}(p)F_{y_0}(p)}{2F_{x_0}(p)},
		\end{align*}
		which gives the osculating $(1,1)$-curve in this case.
		
		If both $b_j\neq 0$ (or both $a_i \neq 0$), or all $a_i,b_j \neq 0$, a linear change of coordinates and straight forward repeated applications of the Euler relations from \cref{euler2} results in the remaining formulas.
\end{proof}
\end{theorem}

\subsection{$(1,1)$-Weierstrass points}
Taking the strategy from \cref{thm:osc11} one step further, we find a criterium for $T_p^{(1,1)}$ to be hyperosculating at $p$, or equivalently, that $p$ a $(1,1)$-Weierstrass point.
\begin{proposition}\label{thm:loccrit11}
	A smooth point $p=(0:1;0:1)$ on a curve $C=\V(F)$ is a $(1,1)$-Weierstrass point if ${F_{x_0}(p),F_{y_0}(p) \neq 0}$ and 
\[\frac{F_{x_0^2y_0}(p)-2k_2\gamma_{00}-2k_4F_{y_0}(p)}{2F_{x_0}(p)} = \frac{F_{x_0y_0^2}(p)-2k_1\gamma_{00}-2k_3F_{x_0}(p)}{2F_{y_0}(p)},\]
	where 
	\[
	\begin{array}{ll}
	k_1=\frac{F_{y_0y_0}(p)}{2F_{y_0}(p)},& \qquad k_2=\frac{F_{x_0x_0}(p)}{2F_{x_0}(p)},\\
	&\\
	k_3=\frac{F_{y_0^3}(p)}{6F_{y_0}(p)},& \qquad k_4=\frac{F_{x_0^3}(p)}{6F_{x_0}(p)},\\
	&\\
	\end{array}\]
	\[\gamma_{00}=F_{x_0y_0}(p)-\frac{F_{y_0y_0}(p)F_{x_0}(p)}{2F_{y_0}(p)}-\frac{F_{x_0x_0}(p)F_{y_0}(p)}{2F_{x_0}(p)}.\]

\begin{proof}
	With notation as in \cref{thm:osc11}, put $P=f-t-k_1ty_0-k_2tx_0-k_3ty_0^2-k_4tx_0^2-k_5tx_0y_0$. Then
	\begin{align*}
	P=H(4)&+\bigl(\tfrac{1}{6}F_{x_0^3}(p)-k_4F_{x_0}(p)\bigr)x_0^3\\
	&+ \bigl(\tfrac{1}{2}F_{x_0^2y_0}(p)-k_2\gamma_{00}-k_4F_{y_0}(p)-k_5F_{x_0}(p)\bigr)x_0^2y_0\\
	&+ \bigl(\tfrac{1}{2}F_{x_0y_0^2}(p)-k_1\gamma_{00}-k_3F_{x_0}(p)-k_5F_{y_0}(p)\bigr)x_0y_0^2\\
	&+\bigl(\tfrac{1}{6}F_{y_0^3}(p)-k_3F_{y_0}(p)\bigr)y_0^3.
	\end{align*}
	
	If the four coefficients in the above expression vanish, then $(T_p^{(1,1)} \cdot C)_p > 3$. This is the case whenever 
	$$\gamma_{00}=F_{x_0y_0}(p)-\frac{F_{y_0y_0}(p)F_{x_0}(p)}{2F_{y_0}(p)}-\frac{F_{x_0x_0}(p)F_{y_0}(p)}{2F_{x_0}(p)},$$
	\[
	\begin{array}{ll}
	k_1=\frac{F_{y_0y_0}(p)}{2F_{y_0}(p)},& \qquad k_2=\frac{F_{x_0x_0}(p)}{2F_{x_0}(p)},\\
	&\\
	k_3=\frac{F_{y_0^3}(p)}{6F_{y_0}(p)},& \qquad k_4=\frac{F_{x_0^3}(p)}{6F_{x_0}(p)},\\
	\end{array}\]
	\[k_5=\frac{F_{x_0^2y_0}(p)-2k_2\gamma_{00}-2k_4F_{y_0}(p)}{2F_{x_0}(p)} = \frac{F_{x_0y_0^2}(p)-2k_1\gamma_{00}-2k_3F_{x_0}(p)}{2F_{y_0}(p)}.\]
	
	\noi Hence, the point $p$ is a $(1,1)$-Weierstrass point whenever such a $k_5$ exists.
\end{proof}
\end{proposition}

\begin{remark} \label{local11hessian}
Consider the expression from \cref{thm:loccrit11} without the evaluation and rewritten as a polynomial, i.e. 
\begin{align*}
H=\;&F_{x_0}^4 (2F_{y_0}F_{y_0^3}-3F_{y_0y_0}^2)\\
&-6F_{x_0}^3F_{y_0}(F_{y_0}F_{x_0y_0^2}-F_{x_0y_0}F_{y_0y_0}) \nonumber \\
&+6F_{x_0}F_{y_0}^3(F_{x_0}F_{x_0^2y_0}-F_{x_0x_0}F_{x_0y_0})\nonumber \\
&-F_{y_0}^4 (2F_{x_0}F_{x_0^3}-3F_{x_0x_0}^2)\nonumber 
\end{align*}
Then in an affine neighbourhood of $(0:1;0:1)$, the curve given by $H(x_0,1,y_0,1)=0$ intersects $C$ at the $(1,1)$-Weierstrass points in the neighbourhood, i.e. it serves as a local $(1,1)$-Hessian. On a side note, notice how in the first and the fourth line the Schwartzian derivative turns up.
\end{remark}

\begin{remark}
	Although \cref{thm:loccrit11} gives a local criteria to determine if a point on a curve is $(1,1)$-Weierstrass and locally acts like a $(1,1)$-Hessian curve, we have not succeeded in finding a global $(1,1)$-Hessian curve. This would be an interesting future task. 
\end{remark}

\subsection{$(1,1)$-Weierstrass point formula}
 In this section we prove a formula for the number of smooth $(1,1)$-Weierstrass points on a curve in $\PP$. For simplicity, we restrict our result to cuspidal curves; a formula for the number of smooth $(1,1)$-Weierstrass points on a curve with general singularities follows by summing over all branches of all singular points. As in \cref{eq:formula1001}, we show the formula by applying results by Ballico and Gatto in \cite{BG97} and Notari in \cite{Not99}, see also \cite{MM17}.
\begin{theorem}[The number of smooth $(1,1)$-Weierstrass points]
    \label{theorem:formula}
Let $C$ be a cuspidal curve of type $(a,b)$ in $\PP$. For a point $p$ let $m_p$ denote the multiplicity and $\delta_p$ the delta invariant. Let $I$ denote the set of cusps that do not have a tangent fiber, and let $\c_p$ denote the intersection multiplicity with an osculating $(1,1)$-curve. Let $J$ denote the set of points that have a tangent fiber, and let $l_p$ denote the intersection multiplicity of the curve and the tangent fiber at $p$. Then the number $W_{(1,1)}$ of smooth $(1,1)$-Weierstrass points on $C$ counted with multiplicity is given by the formula
     \begin{align*}
            W_{(1,1)}=12ab-8a-8b-12\sum_{p \in C} \delta_p-\sum_I (3m_p+\c_p-6)-\sum_J(2m_p+2l_p-6),
    \end{align*}
\begin{proof}
Consider the complete linear system $Q$ of $(1,1)$-curves on $C$. It has dimension $r=r(1,1)=3$, and $\deg Q=a+b$. By \cite[Proposition~3.4]{BG97}, the sum of the Weierstrass weights $w_p(Q)$ of the points on $C$ equals
\begin{align*}
\sum_{p\in C} w_p(Q) = 4(a+b + 3g - 3).
\end{align*}

By \cref{for:genus}, the right hand side can be rewritten as  
\begin{align*}
\sum_{p\in C} w_p(Q) &= 4\bigl(a+b + 3(a-1)(b-1)-3\sum_{p \in C} \delta_p - 3\bigr)\\
&=12ab-8a-8b-12\sum_{p \in C} \delta_p.
\end{align*}

On the other hand, $w_p(Q)$ can be computed as 
\begin{equation}
\label{eq:wp-used}
w_p(Q) = \sum_{i=0}^3 (h_i - i),
\end{equation}
where $h_i$ are distinct integers with the property that $h_i=(C \cdot C_i)_p$ for curves $C_0,\dots,C_3$ of type $(1,1)$.

To compute $h_i$ we make use of the Puiseux parametrization of $C$ at a point $p$ \cite[Cor.~7.7, p.~135]{Fis01} and the standard basis of curves of type $(1,1)$, \[x_0y_0, x_0y_1,x_1y_0,x_1y_1.\] 

In $\PP$ there are two cases to consider. If a fiber is tangent to $C$ at $p$, the curve can locally be given by the parametrization
\begin{equation*}
\label{eq:puiseux_paramfibtan}
(t^m:1; a_lt^l +\dots : 1),
\end{equation*}
where $m=m_p$ and $l=l_p$ are as above, $a_l\neq0$, and $"\dots"$ denotes higher order terms in $t$. Substituting the parametrization into the standard basis, we find \[a_lt^{m+l}+\dots, t^m,a_lt^l+\dots,1,\] and subsequently \[h_0=0,h_1=m,h_2=l,h_3=m+l.\] By \cref{eq:wp-used} \[w_p(Q)=2m+2l-6.\]

If a fiber is not tangent to $C$ at $p$, the curve can locally be given by the parametrization
\begin{equation*}
\label{eq:puiseux_param}
(t^m :1; a_mt^m+a_lt^{l}\dots : 1),
\end{equation*}
where $m<l$, and $a_m,a_l \neq 0$. Substituting the parametrization into the standard basis, we get \[a_mt^{2m}+a_lt^{m+l}+\dots, t^m,a_mt^m+a_lt^{l}+\dots,1.\] In this case $h_0=0$, $h_1=m$ and $h_2=2m$. Since $a_l \neq 0$, by changing the basis we may construct a $(1,1)$-curve $C_3$ such that $h_3=\c$ for an integer $\c \neq m,2m$, referred to as an osculating $(1,1)$-curve at $p$. This curve coincides with $T_p^{(1,1)}$ if $p$ is smooth. Note that $l$ can be interpreted as the intersection multiplicity of the local tangent line to $C$ at $p$ in an affine neighbourhood of $p$. Moreover, $\c=l$ if $l \neq 2m$. By \cref{eq:wp-used} \[w_p(Q)=3m+\c-6.\]

For a smooth point that does not have a tangent fiber we have that $m=1$, and thus \[w_p(Q)=\c-3,\] which equals its multiplicity as a $(1,1)$-Weierstrass point. 

Summing up and reorganizing terms, we reach the desired expression.
\end{proof}
\end{theorem}

Note that the integers $l_p$ and $c_p$ can be estimated by, and frequently determined by, the multiplicity sequence of the cusp $p$, see \cite[p.~503]{BK86}.
\begin{lemma} Let $\ol{m}_p=[m,m_1,m_2,\dots]$ denote the multiplicity sequence of a cusp $p$. If a fiber is tangent to $C$ at $p$, then depending on the direction either \[a \geq l_p=k\cdot m+m_k\] or \[b \geq l_p=k\cdot m+m_k,\] for some $k \geq 1$ with $m=m_1=\dots=m_{k-1}$. 
	
If no fiber is tangent to $C$ at $p$, then the value of $c_p$ is bounded by \[a+b \geq c_p=k\cdot m+m_k,\] for some $k \geq 1$ with $m=m_1=\dots=m_{k-1}$.
	\begin{proof}
		The inequalities on the left hold by taking the intersection of $C$ and curves of respective types. The equalities on the right and the restrictions on $k$ and the multiplicity sequence follow from \cite[Proposition 1.4]{FZ96} since $(1,1)$-curves are non-singular.
	\end{proof}
	
\end{lemma}

\begin{remark}
	\label{rem:11hes} 	
	The formula in \cref{theorem:formula} indicates that a $(1,1)$-Hessian curve, if it exists, should be of type $(6a-8,6b-8)$. Indeed, intersecting a curve of type $(a,b)$ with a curve of type $(6a-8,6b-8)$ gives in total $12ab-8a-8b$. In this case, the formula from \cref{theorem:formula} has, as in \cref{conj:1001}, an interpretation as an application of Bézout's theorem for $\PP$, where the sums can be reorganized to represent local intersection multiplicities.
	\end{remark}

\section{The special case of rational curves}\label{sec:rational}
In this section we study rational curves in $\PP$. For any rational curve $C$, the geometric genus $g=0$, and the curve can be given by a parametrization $\Phi(s,t)$,
\[
\Phi(s,t) = (\varphi_0(s,t) : \varphi_1(s,t) ; \psi_0(s,t) : \psi_1(s,t)) \; \text{for } (s:t) \in \mathbb{P}^1,
\]
where $\deg \varphi_0 = b=\deg \varphi_1$ and $\deg \psi_0 = a = \deg \psi_1$.
\subsection{Osculating curves}
We first observe that the fibers through a point can be expressed by the parametrization, which follows immediately from \cref{osc1001curve}.
\begin{proposition}
Let $C$ be a rational curve in $\PP$ given by $\Phi(s,t)$, and let \[\omega_{(1,0)}(s,t)=\varphi_1(s,t)x_0-\varphi_0(s,t)x_1,\] and \[\omega_{(0,1)}(s,t)=\psi_1(s,t)y_0-\psi_0(s,t)y_1.\] Then for a smooth point $p=\Phi(s_0,t_0)$, the polynomial $\omega_{(1,0)}(s_0,t_0)$ is the defining polynomial of the osculating $(1,0)$-curve to $C$ at $p$. Symmetrically, the polynomial $\omega_{(0,1)}(s_0,t_0)$ is the defining polynomial of the osculating $(0,1)$-curve to $C$ at $p$.
\end{proposition}

We may find a similar expression for the osculating $(1,1)$-curve at any given smooth point of a curve. 
\begin{proposition}
Let $C$ be a rational curve in $\PP$ given by $\Phi(s,t)$, and let $\omega_{(1,1)}(s,t)$ denote the determinant
\begin{align*}
\omega_{(1,1)}(s,t) = \begin{vmatrix}
x_0y_0 & x_0y_1 & x_1y_0 & x_1y_1 \\[1ex]
\frac{\di^2}{\di s^2}(\varphi_0\psi_0) & 
\frac{\di^2}{\di s^2}(\varphi_0\psi_1) & 
\frac{\di^2}{\di s^2}(\varphi_1\psi_0) & 
\frac{\di^2}{\di s^2}(\varphi_1\psi_1)  \\[1ex]
\frac{\di^2}{\di s \di t}(\varphi_0\psi_0) & 
\frac{\di^2}{\di s\di t}(\varphi_0\psi_1) & 
\frac{\di^2}{\di s\di t}(\varphi_1\psi_0) & 
\frac{\di^2}{\di s\di t}(\varphi_1\psi_1)  \\[1ex]
\frac{\di^2}{\di t^2}(\varphi_0\psi_0) & 
\frac{\di^2}{\di t^2}(\varphi_0\psi_1) & 
\frac{\di^2}{\di t^2}(\varphi_1\psi_0) & 
\frac{\di^2}{\di t^2}(\varphi_1\psi_1)  \\
\end{vmatrix}.
\end{align*}
Then for a smooth point $p=\Phi(s_0,t_0)$, the determinant $\omega_{(1,1)}(s_0,t_0)$ is the defining polynomial of the osculating $(1,1)$-curve to $C$ at $p$.
\begin{proof}
	Let $s(C)$ denote the image of $C$ under the Segre-embedding $s \colon \PP \longrightarrow \mathbb P^3$, given by the parametrization 
	\begin{align*}
s(C) = (\varphi_0\psi_0: \varphi_0\psi_1 : \varphi_1\psi_0 : \varphi_1\psi_1).
	\end{align*}
	Consider the determinant $\tilde{\omega}(s,t)$, where in the first row of $\omega(s,t)$ the standard basis of $(1,1)$-curves is substituted with the coordinates of $\P^3$. Then for a smooth point $s(C)(s_0,t_0)$ in $\P^3$, the determinant $\tilde{\omega}(s_0,t_0)$ provides a linear polynomial that defines the unique osculating plane to $s(C)$ in $\P^3$. This plane corresponds to the osculating $(1,1)$-curve to $C$ at $p=\Phi(s_0,t_0)$ in $\PP$, with defining polynomial $\omega(s_0,t_0)$.
\end{proof}

\end{proposition}

\subsection{Weierstrass points and weight}
We now locate Weierstrass points on rational curves, and compute their Weierstrass weight, directly from the parametrization.

\begin{theorem}\label{rational1001hessian}
	Let $C$ be a rational curve in $\PP$ given by $\Phi(s,t)$, and let $\xi_{(1,0)}(s,t)$ denote the determinant 
		\begin{align*}
	\xi_{(1,0)}(s,t) = \begin{vmatrix}
	\frac{\di}{\di s}(\varphi_0) &
	\frac{\di}{\di s}(\varphi_1) \\[1ex]
	\frac{\di}{\di t}(\varphi_0) &
\frac{\di}{\di t}(\varphi_1)
	\end{vmatrix}.
	\end{align*}
	Moreover, let $(s_i:t_i)$ denote the distinct zeros of $\xi_{(1,0)}(s,t)$, $i \leq 2b-2$. Then the points $p_i=\Phi(s_i,t_i)$ are the $(1,0)$-Weierstrass points on $C$, and the Weierstrass weight $w_{p_i}(1,0)$, is equal to the order of the zero of $\xi_{(1,0)}(s,t)$ corresponding to $(s_i:t_i)$. A symmetrical result holds for $(0,1)$-Weierstrass points, then with 
	\begin{align*}
	\xi_{(0,1)}(s,t) = \begin{vmatrix}
	\frac{\di}{\di s}(\psi_0) &
	\frac{\di}{\di s}(\psi_1) \\[1ex]
	\frac{\di}{\di t}(\psi_0) &
	\frac{\di}{\di t}(\psi_1)
	\end{vmatrix}.
	\end{align*}
	\begin{proof}
The polynomial of a $(1,0)$-curve through a point $p$ can be written as \[\varphi_1x_0-\varphi_0x_1.\] This polynomial has a multiple zero in $(s_i:t_i)$ whenever $\xi_{(1,0)}(s_i,t_i)$ vanishes. The order of the zero in the polynomial is equal to the intersection multiplicity of the curve and the fiber at the point, and the order of the corresponding zero of $\xi_{(1,0)}(s,t)$ is one less, hence equal to the $(1,0)$-Weierstrass weight from \cref{eq:formula1001}.
	\end{proof}
\end{theorem}

\begin{remark}
Since $\deg \xi_{(1,0)}(s,t)=2(b-1)$ and $\deg \xi_{(0,1)}(s,t)=2(a-1)$, \cref{rational1001hessian} gives a direct proof of \cref{eq:formula1001} in the rational case.
\end{remark}

\begin{theorem}\label{rational11hessian}
	Let $C$ be a rational curve of type $(a,b)$ in $\PP$ given by $\Phi$, and let $\xi_{(1,1)}(s,t)$ denote the determinant
	\begin{align*}
	\xi_{(1,1)}(s,t) = \begin{vmatrix}
	\frac{\di^3}{\di s^3}(\varphi_0\psi_0) & 
	\frac{\di^3}{\di s^3}(\varphi_0\psi_1) & 
	\frac{\di^3}{\di s^3}(\varphi_1\psi_0) & 
	\frac{\di^3}{\di s^3}(\varphi_1\psi_1) \\[1ex]
	\frac{\di^3}{\di s^2\di t}(\varphi_0\psi_0) &
	\frac{\di^3}{\di s^2\di t}(\varphi_0\psi_1) &
	\frac{\di^3}{\di s^2\di t}(\varphi_1\psi_0) &
	\frac{\di^3}{\di s^2\di t}(\varphi_1\psi_1) \\[1ex]
	\frac{\di^3}{\di s\di t^2}(\varphi_0\psi_0) &
	\frac{\di^3}{\di s\di t^2}(\varphi_0\psi_1) &
	\frac{\di^3}{\di s\di t^2}(\varphi_1\psi_0) &
	\frac{\di^3}{\di s\di t^2}(\varphi_1\psi_1) \\[1ex]
	\frac{\di^3}{\di t^3}(\varphi_0\psi_0) &
	\frac{\di^3}{\di t^3}(\varphi_0\psi_1) & 
	\frac{\di^3}{\di t^3}(\varphi_1\psi_0) & 
	\frac{\di^3}{\di t^3}(\varphi_1\psi_1)
	\end{vmatrix}.
	\end{align*}
	Moreover, let $(s_i:t_i)$ denote the distinct zeros of $\xi_{(1,1)}(s,t)$, $i \leq 4(a+b-3)$. Then the points $p_i=\Phi(s_i,t_i)$ are the $(1,1)$-Weierstrass points on $C$, and the Weierstrass weight $w_{p_i}(1,1)$, is equal to the order of the zero of $\xi_{(1,1)}(s,t)$ corresponding to $(s_i:t_i)$.
\begin{proof}
	The zeros of $\xi(s,t)$ correspond to points of the curve $s(C)$ with a hyperosculating plane in $\P^3$, hence points on $C$ with a hyperosculating $(1,1)$-curve in $\PP$. Thus they correspond to the $(1,1)$-Weierstrass points of $C$, and there are at most $\ord \xi(s,t) \leq 4(a+b-3)$ such points. 

It remains to show that the order of the zeros equals the Weierstrass weight of the point. 

Assume first that $p$ on $C$ is unibranched, i.e. there is a unique pair $(s:t)$ such that $p=\Phi(s,t)$. Using that a unibranched point either has a tangent fiber or not, we compute the Weierstrass weight directly using the corresponding Puiseux parametrizations. 
	
	If a fiber is tangent to $p=\Phi(1,0)$, consider the local parametrization \[(t^m:1;a_lt^l +\dots:1),\] with $m<l$. Then \[\xi_{(1,1)}(1,t)=(m-l)(m+l)m^2l^2a_l^2t^{2m+2l-6}+\dots,\] where the coefficient is non-zero, so by comparing with the proof of \cref{theorem:formula} \[\ord_t\xi=2m+2l-6=w_p(1,1).\]
	
	If a fiber is not tangent to $p=\Phi(1,0)$, then first write the local parametrization \[(t^m:1;a_mt^m+a_lt^l+\dots:1),\] where $a_m,a_l \neq 0$ and $m<l$, and $l$ can be interpreted as the intersection multiplicity of $C$ and the classical tangent line in the affine neighbourhood. Then \[\xi_{(1,1)}(1,t)=2m^3la_ma_lt^{3m+l-6}(l^2-3lm+2m^2)+\dots,\] hence \[\ord_t\xi=3m+l-6 \text{ if } l\neq 2m.\] Moreover, by examining the standard basis of $(1,1)$-curves, there exists a $(1,1)$-curve that intersects $C$ at $p$ with intersection multiplicity $l$. 
	
	If $l=2m$, then since $p$ is unibranced, we may write the local parametrization \[(t^m:1;a_mt^m+a_{2m}t^{2m}+a_ct^c+\dots:1),\] for an integer $c$ such that $2m<c$ and $a_c \neq 0$. Then similarly, \[\xi(1,t)=2m^3ca_ma_ct^{3m+c-6}(c^2-3cm+2m^2)+\dots,\] and since $c>2m>m$ the coefficient is non-zero and \[\ord_t\xi=3m+c-6 \text{ if } l=2m.\] Again by examining the standard basis of $(1,1)$-curves, there exists a $(1,1)$-curve that intersects $C$ at $p$ with intersection multiplicity $c$.
		
	In either case, there exists an integer $c$ such that it is the intersection multiplicity of a $(1,1)$-curve and $C$ at $p$, $c=l$ if $l \neq 2m$, and comparing with the proof of \cref{theorem:formula} \[\ord_t\xi=3m+c-6=w_p(1,1).\]
	
	Now, if $p$ is a multiple point, we compute the Weierstrass weight for each value of $(s_i:t_i)$ such that $p=\Phi(s_i,t_i)$ separately, and sum up to reach $w_p(1,1)$.
	
\end{proof}	
\end{theorem}

\begin{remark}
	Since $\deg \xi_{(1,1)}(s,t)=4(a+b)-3$, \cref{rational11hessian} gives a direct proof of \cref{theorem:formula} in the rational case.
\end{remark}

\section{Applications}\label{sec:applications}
In this section we apply the results in this article to a curve $C$ in $\PP$. In the first example we take a closer look at a rational curve.
\begin{example}
Let $C$ be the curve in $\PP$ of type $(3,2)$ given by the parametrization \[\Phi=(-st+t^2:s^2;t^3:s^3),\] with defining polynomial \[F=x_0^3y_1^2+3x_0x_1^2y_0y_1-x_1^3y_0^2+x_1^3y_0y_1.\] 

See \cref{tab:exa} for an overview of the singular points and the Weierstrass points for this curve. In the table, $\zeta_x=\tfrac{5}{18}\zeta_y-\tfrac{2}{9}$, $\zeta_y=\tfrac{136}{125}\pm \tfrac{54}{125}\sqrt{6}$, and $\zeta_t = \tfrac{4}{5}\pm \tfrac{1}{5}\sqrt{6}$.

		\begin{table}[ht]
	\centering
	\resizebox{\textwidth}{!}{
		\begin{tabular}{ccccccccc}
			
			Point $p$ & $(s:t)$& Point type & $m_p$ & $\delta_p$ & $\#$ branches &$w_p(1,0)$&$w_p(0,1)$ &$w_p(1,1)$\\
			\midrule
	$(1:0;1:0)$&$(0:1)$ & Cusp & $2$ & $1$ & $1$ & $1$ & $2$ & $4$\\
	$(-1:1;-1:1)$&$(2:1\pm\sqrt{3}i)$ & Node & $2$ & $1$ & $2$ & $0$ & $0$ & $0$\\
	$(0:1;0:1)$&$(1:0)$ & Higher order $(0,1)$-WP & $1$ & $0$ & $1$ & $0$ & $2$ & $2$\\
	$(-1:4;1:8)$&$(2:1)$ & Ordinary $(1,0)$-WP & $1$ & $0$ & $1$ & $1$ & $0$ & $0$\\
	$(\zeta_x:1;\zeta_y:1)$&$(1:\zeta_t)$ & Two smooth $(1,1)$-WPs & $1$ & $0$ & $1$ & $0$ & $0$ & $1$\\
		\end{tabular}
	}
	\caption{An overview of the curve $C$.}
	\label{tab:exa}
\end{table}

The Weierstrass points and the Weierstrass weights in \cref{tab:exa} can be found by computing the determinants from \cref{sec:rational}. The Weierstrass points can be read off the factorization, and the Weierstrass weights are equal to the multplicities of the corresponding factors;
\begin{align*}
\xi_{(1,0)}&=2(s-2t)s=0,\\
\xi_{(0,1)}&=-9s^2t^2=0,\\
\xi_{(1,1)}&=-77760s^4t^2(2s^2-8st+5t^2)=0.
\end{align*} 

For any given point $p=\Phi(s,t)$, the corresponding osculating $(1,1)$-curve is given by $\omega_{(1,1)}(s,t)$:
\begin{equation*}
\resizebox{.98 \textwidth}{!} 
{$
s^2t\bigl[s^5(s-t)x_0y_0+s^2t^3(2s-5t)x_0y_1+s^3t(2s^2-6st+5t^2)x_1y_0+t^4(s^2-3st+t^2)x_1y_1\bigr].
$}
\end{equation*}
Observe that the polynomial vanishes for $(s:t)=(0:1)$ and $(1:0)$, which correspond to the cusp and the higher order $(0,1)$-Weierstrass point of $C$, where the osculating $(1,1)$-curve in each case degenerates to the product of the fibers through the point. 

In particular, note that for the $(1,1)$-Weierstrass point given by the parameters $(1:\tfrac{4}{5}+ \tfrac{1}{5}\sqrt{6})$, 
\begin{align*}
\omega_{(1,1)}(1,\tfrac{4}{5}+ \tfrac{1}{5}\sqrt{6})=&\;3125x_0y_0\\
&+(21000\sqrt{6}+51500)x_0y_1\\
&-(7500x_1y_0\sqrt{6}+17500)x_1y_0\\
&+(7344\sqrt{6}+17996)x_1y_1,
\end{align*}
which can also be found using \cref{thm:osc11}. 

The number of Weierstrass points can be verified by \cref{eq:formula1001} and \cref{theorem:formula} after a computation of the singular points and connected invariants.
\begin{align*}
	W_{(1,0)}&=2(b+g-1)-\sum(l_{p_j}-1)\\
&=2(2+0-1)-(2-1)\\
&=1.\\
	W_{(0,1)}&=2(a+g-1)-\sum (l_{p_j}-1)\\
			 &=2(3+0-1)-(3-1)\\
			 &=2.\\
W_{(1,1)}&=12ab-8a-8b-12\sum_{p \in C} \delta_p-\sum_I (3m_p+\c_p-6)-\sum_J(2m_p+2l_p-6)\\
&=12\cdot 3\cdot 2-8\cdot 3-8\cdot 2-12 \cdot 2-(2\cdot 2+2 \cdot 3 -6)-(2 \cdot 1+2 \cdot 3-6)\\
&=2.
\end{align*} 
\end{example}

\begin{remark}
We have by direct computation verified that the conjectured results from \cref{conj:1001} and \cref{rem:11hes} hold in this case, i.e. that for each point the local intersection multiplicities of the curve and the corresponding Hessian in each case equals the sums of the repsective Weierstrass weights and multiple of the delta invariants. 	
\end{remark}

In the next example we generalize the method from \cref{sec:rational} and show that a special class of curves does not have smooth $(\a,\b)$-Weierstrass points.
\begin{example}
Let $C$ be the curve of type $(a,b)$, where $2 \leq a<b$ and $\gcd(a,b)=1$, given by the parametrization \[\Phi=(s^b:t^b;s^a:t^a),\] with defining polynomial \[F=x_0^ay_1^b-x_1^ay_0^b.\] The curve has two cusps $p_1=\Phi(0,1)$ and $p_2=\Phi(1,0)$ of the same analytic type and with multiplicities $m_{p_1}=a=m_{p_2}$.

Let $(\a,\b)$, subject to the restrictions $0 \leq \a<b$ and $0\leq \b<a$, be the complete linear system of curves of type $(\a,\b)$ on $C$. Then substituting the parametrization in the standard basis for $H^0(\PP,\mathcal{O}_{\PP}(\a,\b))$, we get a list of $r+1$ different monomials of degree $b\a+a\b$ in $s$ and $t$. Computing the corresponding $(\a,\b)$-Wronskian leads to \[\xi_{(\a,\b)}(s,t)=Ks^nt^n,\]
where $K \in \N$ is a function of $a,b,\a,\b$, and $n=\frac{(\a+1)(\b+1)\bigl(\b(a-1)+\a(b-1)-\a\b\bigr)}{2}$. Indeed the order of $s$ and $t$ in $\xi_{(\a,\b)}(s,t)$ must be equal, and the total degree of $\xi_{(\a,\b)}(s,t)$ is equal to ${(r+1)(b\a+a\b-r)}$, and substituting for $r$, we see that each order must be $n$. 

Thus the cusps are the only $(\a,\b)$-Weierstrass points on $C$, their $(\a,\b)$-Weierstrass weight is $n$, and there are no smooth $(\a,\b)$-Weierstrass points on $C$.
\end{example}	

\subsection*{Acknowledgements}
	This article is based on methods from the master's thesis of the first author supervised by the second author and Geir Ellingsrud \cite{Mau17}, and previously unpublished results from the PhD thesis of the second author \cite{Moe13}. We would like to thank Hubert Flenner for requesting a Plücker-like formula for $\PP$, and Kristian Ranestad for his support. Moreover, we would like to thank Head of the Science Library at the University of Oslo, Live Rasmussen, for the opportunity to write this paper.

\bibliographystyle{amsalpha2}
\bibliography{bibPP}{}

\end{document}